\documentstyle[righttag,amssymb,intlim,amscd,bezier]{amsart}

\subjclass{53C12, 70H05}
\keywords{Foliation, vertical cohomology, Hamiltonian systems,
    isoenergetic surfaces, Liouville torus}

\newcommand{\nc}{\newcommand}
\nc{\al}{\alpha}
\nc{\bt}{\beta}
\nc{\gm}{\gamma}
\nc{\dl}{\delta}
\nc{\Lb}{\Lambda}
\nc{\vf}{\varphi}
\nc{\om}{\omega}
\nc{\tht}{\theta}

\nc{\bR}{\Bbb{R}}
\nc{\bT}{\Bbb{T}}
\nc{\bZ}{\Bbb{Z}}
\nc{\cF}{\cal{F}}
\nc{\cP}{\cal{P}}
\nc{\cT}{\cal{T}}

\nc{\pa}{\partial}
\nc{\sbs}{\subset}
\nc{\wh}{\widehat}
\nc{\wt}{\widetilde}
\nc{\ol}{\overline}
\nc{\os}{\overset}
\nc{\us}{\underset}
\nc{\str}{\stackrel}
\nc{\map}{\longmapsto}
\nc{\hook}{\hookrightarrow}

\nc{\sgrad}{\operatorname{sgrad}}
\nc{\grad}{\operatorname{grad}}
\nc{\Inv}{\operatorname{Inv}}
\nc{\End}{\operatorname{End}}

\newtheorem{th}{\ \ \ Theorem}[section]
\newtheorem{prop}[th]{\ \ \ Proposition}
\newtheorem{cor}[th]{\ \ \ Corollary}
\newtheorem{lem}[th]{\ \ \ Lemma}

\theoremstyle{definition}
\newtheorem{defi}[th]{\ \ \ Definition}

\theoremstyle{remark}
\newtheorem{rem}[th]{\ \ \ Remark}

\numberwithin{equation}{section}

\begin{document}

\title[Vertical Cohomologies and Their Application]
    {Vertical Cohomologies and Their Application to
        Completely  Integrable Hamiltonian  Systems}

\author{Z. Tevdoradze}

\maketitle
\bigskip

\hfill {\it Dedicated to the memory of Roin Nadiradze}
\vskip-2cm

\begin{abstract}
    Some functorial and topological properties of vertical
    cohomologies and their application to completely integrable
    Hamiltonian systems are studied.
\end{abstract}

%   \S 1.
\section{Introduction }

    If on a smooth Riemannian manifold $M^n$ we have a distribution $V$ of
dimension $k$, which is actually a smooth section of the Grassman
fiber bundle $G_k(TM^n)\to M^n$ adjoint to the tangent fibration
$TM^n$ to the manifold $M^n$, then by means of the Riemannian
metric we obtain $TM^n=V\oplus N$, where $N$ is a normal fiber
bundle to $V$. Let $P:TM^n\to V\sbs TM^n$ be a natural projection.
The operator $P$ defines the mapping $P^*:\Lb^*(M^n)\to
\Lb^*(M^n)$ ($(\Lb^*(M^n),d^*)$ is the de Rham differential
complex) by the formula
$(P^*\al)(X_1,\dots,X_q)=\al(PX_1,\dots,PX_q)$, where $\al\in
\Lb^q(M^n)$ and $X_1,\dots,X_q\in S(M^n)$ are the smooth vector
fields on $M^n$.

    Denote by $\Lb^*_V(M^n)$ all fixed points of the operator $P^*$.
In what follows we shall consider the case, where $V$ is
integrable, i. e., where $M^n$ is partitioned into leaves and the
tangent space to the leaf that passes through the point  $x\in
M^n$ is $V_x$. Then the pair $(\Lb^*_V(M^n),d_V^*)$ forms a
differential complex with the differential $d_V^*=P^*\circ d^*$.

    The cohomologies  of the complex $(\Lb_V^*(M^n),d_V^*)$ are
called vertical and denoted by $H_V^*(M^n)$ ([1]). It has turned
out that these cohomologies coincide with those of the classical
BRST operator ([2], [3]).

    In \S 2 the vertical cohomologies are defined without fixing
the metric on $M^n$, and some of their functorial properties are
studied. The FOL category  of smooth foliations and leaf-to-leaf
transforming mappings is introduced, and a natural transformation
of the de Rham functor $\Lb^*$ to the functor $\Lb^*_{\cF}$ is
construced (Proposition 2.2). The notion of leaf-to-leaf
transforming homotopic mappings is introduced, and the homotopy
axiom for vertical cohomologies is proved (Theorem 2.5). The
notion of  a relative group of vertical cohomologies is introduced
by analogy with de Rham's theory, and the long exact cohomologic
sequences (2.6) and (2.7) are derived. Moreover, for a
leaf-to-leaf transforming mapping $f:(M^n,{\cF}_1)\to
(N^m,{\cF}_2)$, the cohomology groups $H^*(f)$ are constructed and
proved (Theorem 2.8) to be isomorphic for leaf-to-leaf
transforming homotopic mappings. Finally, a double complex
$(K^{**},D^*)$ is constructed for the countable covering
$U=\{u_\al\}_{\al\in A}$ of the manifold $M^n$. It is shown that
the cohomologies of $(K^{**},D^*)$ are isomorphic to the vertical
cohomologies. A combinatorial definition of vertical comologies in
the \v{C}ech sense (Theorems 2.10 and 2.12) is also given.

    In \S 3 some of the main facts from the topological theory of
integrable  Hamiltonian systems  ([5], [6]) are presented. Using
the notion of vertical cohomologies, the groups corresponding to
nonresonance Hamiltonian systems are constructed (Theorems 3.3 and
3.4).

    In \S 4 the case of a spherical pendulum is considered as an example.

%   \S 2.
\section{Vertical Cohomologies}

%   2.1.
\subsection{Definition of Vertical Cohomologies.}
 Let $M^n$ be
a smooth $n$-dimensional manifold, and $V$ be a $k$-dimensional
involutive distribution on $M^n$ whose foliation is denoted by
$\cF$.  The bundle of exterior $p$-forms on $V$ is denoted by
$A^p(V)$, and the  set of  smooth sections of the bundle $A^p(V)$
by $\Lb^p_{\cF}(M^n)$. Then $\Lb_{\cF}^p(M^n)$ is a module over
the algebra of infinitely differentiable functions $C^\infty(M^n)$
on $M^n$.

    Let $\al\in \Lb_{\cF}^p(M^n)$, and let $X_1,\dots,X_p$ be the smooth
vector fields on $M^n$ which are tangent to the leaves of the
foliation $\cF$, i.e., they are the smooth sections of the bundle
$V\str{p_{{}_{\cF}}}{\to} M^n$. Then the mapping  defined by the
formula
%   (2.1)
\begin{equation}
    \al(X_1,\dots,X_p):x\str{j}{\map} \al(X_1(x),\dots,X_p(x)),
\end{equation}
which on the module of  sections $S(V)$ of $V$ assigns an exterior
$p$-form to an element $\al\in \Lb_{\cF}^p(M^n)$, is an
isomorphism.

    Let us now define the operator
$$  d_{\cF}^p:\Lb_{\cF}^p(M^n)\to \Lb_{\cF}^{p+1}(M^n)  $$
by the relation
\begin{gather*}
    (d_{\cF}^p\al)(X_1,\dots,X_{p+1}) =\sum_{i=1}^{p+1}
        (-1)^{i-1} X_i\al(X_1,\dots,\wh{X}_i,X_{i+1},
        \dots,X_{p+1})+  \\
    +\sum_{i<j} (-1)^{i+j} \al([X_i,X_j],X_1,\dots,\wh{X}_i,\dots,
        \wh{X}_j,\dots,X_{p+1}),
\end{gather*}
where $\al\in \Lb_{\cF}^p(M^n)$, $X_i\in S(V)$, $i=\ol{1,p+1}$.

        It is easy to verify that the embedding
$i:S(V)\hook S(TM^n)$ induces the projection $i_q^*:A^q(S(TM))\to
A^q(S(V))$. Indeed, let $\al\in A^q(S(V))$. Define $\ol{\al}\in
A^q(S(TM))$ so that $i_q^*\ol{\al}=\al$. Any Riemannian metric
defines a smooth section $P:M^n\to \End (TM^n)$ of the bundle of
homomorphisms $\End(TM^n)$ of the tangent bundle $TM^n$, where
$P(x):T_xM^n\to V_x\sbs T_xM^n$ is the orthogonal  projection
onto  a leaf of the foliation $\cF$ which passes through the point
$x\in M^n$. The element $\ol{\al}$ can be defined by the formula
$$  \ol{\al}(X_1,\dots,X_q)=\al(PX_1,\dots,PX_q),\quad
        X_i\in S(TM^n),\;\;\;i=\ol{1,q}.    $$

        The following diagram is commutative:
$$  \begin{CD}
        A^q(S(TM^n)) @>{i_q^*}>> A_{\cF}^q(M^n) \\
        @VV{d^q}V @VV{d_{\cF}^q}V \\
        A^{q+1}(S(TM^n)) @>{i_{q+1}^*}>> A_{\cF}^{q+1}(M^n)
    \end{CD}\;,     $$
where by $A_{\cF}^*(M^n)$ is denoted $A^*(S(V))$. Indeed,
\begin{gather*}
    (d_{\cF}^qi_q^*\al)(X_1,\dots,X_{q+1}) =\sum_{j=1}^{q+1}
        (-1)^{j-1} (X_j)(i_q^*\al)(X_1,\dots,\wh{X}_j,
        \dots,X_{q+1})+  \\
    +\sum_{j<t} (-1)^{j+t} (i_q^*\al)([X_j,X_t],X_1,\dots,\wh{X}_j,\dots,
        \wh{X}_t,\dots,X_{q+1})= \\
    =(d^q\al)(X_1,\dots,X_{q+1})=(i_{q+1}^* d^q\al)(X_1,\dots,X_{q+1}).
\end{gather*}

      The commutativity $d_{\cF}^q\circ i_q^*=i_{q+1}^*\circ d^q$
implies $d_{\cF}^q\circ d_{\cF}^{q-1}=0$. Thus $i^*$ is a cochain
mapping between the differential complexes $(A^*(S(TM^n),d^*)$ and
$(A_{\cF}^*(M^n),d_{\cF}^*)$.

%        DEFINITION 2.1.
\begin{defi}
Cohomology groups of the complex $(A_{\cF}^*(M),d_{\cF}^*)$ are
called vertical cohomologies of the foliation $(M^n,\cF)$, and we
denote them by $H_{\cF}^*(M^n)$.
\end{defi}

        It is easy to verify that $d_{\cF}^*$ is an antidifferentiation
of order 1, i.e., if $\al\in A_{\cF}^q(M^n)$, $\bt\in
A_{\cF}^l(M^n)$, we have
$$  d_{\cF}^{q+l}(\al\wedge \bt)=(d_{\cF}^q\al)\wedge \bt+
        (-1)^q\al\wedge (d_{\cF}^l\bt),     $$
where $\wedge$ is an exterior product. This implies that the
homomorphism $i^*$ induces a homomorphism  between the de Rham
cohomology algebra and the cohomology algebra $H_{\cF}^*(M^n)$. We
denote this homomorphism by the same symbol $i^*$. Since we
already know that the cochain mapping $i^*$ is an epimorphism, we
have a short exact sequence of cochain differential complexes
%   (2.2)
\begin{equation}
    0 @>>> Z_{\cF} @>{j^*}>> A^* (S(TM^n)) @>{i^*}>>
        A_{\cF}^*(M^n) @>>> 0,
\end{equation}where by $(Z_{\cF}^*,d^*)$ is denoted the kernel of the mapping
$i^*$. Sequence (2.2) induces a long exact sequence of cohomology
groups
%   (2.3)
\begin{equation}
    0 @>>> \wt{H}_{\cF}^0(M^n) @>{j^*}>> H^0(M^n) @>{i^*}>>
        H_{\cF}^0(M^n) @>{\dl}>> \wt{H}_{\cF}^1(M^n)
        @>{j^*}>> \cdots \;.
\end{equation}
Here by $\wt{H}_{\cF}^*(M^n\!)$ are denoted the cohomology groups
of the complex $(Z_{\cF}^*,d)$. Since $H_{\cF}^{k+m}(M^n)=0$,
$m>0$, from (2.3) we get $\wt{H}_{\cF}^{m+k}(M^n) \approx
H^{m+k}(M^n)$, $m>1$.

        Denote by FOL a category whose objects are smooth foliations
$(M^n,\cF)$, and morphisms from $(M_1^n,{\cF}_1)$ to
$(M_2^m,{\cF}_2)$ are leaf-to-leaf transforming mappings, i.e.,
smooth mappings $h$ from $M_1^n$ to $M_2^n$ which preserve the
foliation structure by perform a leaf-to-leaf transformation.

        If $h$ is the leaf-to-leaf transforming  mapping between the foliations
\linebreak  $(M_1^n,{\cF}_1)$ and $(M_2^n,{\cF}_2)$, then $h$
defines the morphism $h^*$ between the    \linebreak
$C^\infty(M_2^m)$-module $\Lb_{{\cF}_2}^*(M_2^m)$ and the
 $C^\infty(M_1^n)$-module $\Lb_{\cF_1}^*(M_1^n)$
by the formula
%   (2.4)
\begin{equation}
    (h^*\al)(X_1,\dots,X_p)(x)=\al(h^\top X_1(x),\dots,
        h^\top X_p(x)),
\end{equation}
where $\al\in \Lb_{{\cF}_2}^p(M_2)$ and $X_1,\dots,X_p\in S(V_1)$,
$V_1$ is the distribution associated with ${\cF}_1$, and $h^\top$
is the tangent mapping to $h$. We have the commutative diagram
$$  \begin{CD}
        \Lb^*(M_2^m) @>{h^*}>> \Lb^*(M_1^n) \\
        @VV{i_2^*}V @VV{i_1^*}V \\
        \Lb_{{\cF}_2}^*(M_2^m) @>{h^*}>> \Lb_{{\cF}_1}^*(M_1^n)
    \end{CD}\;.     $$

        A direct calculation shows that $h^*$ is a cochain mapping, i.e.,
$d_{{\cF}_1}^*\circ h^*=h^{*+1}\circ d_{{\cF}_2}^*$. Hence we have

%        PROPOSITION 2.2.
\begin{prop}
The mapping $i^*$ is a natural transformation of the de Rham
functor $\Lb^*$ to the functor $\Lb_{\cF}^*$, where $\Lb^*$ and
$\Lb_{\cF}^*$ are the contravariant functors from the \rom{FOL}
category to the category of differential graded algebras and their
homomorphisms.
\end{prop}

%        2.2.
\subsection{A Homotopy Axiom for Vertical Cohomologies.}

        Let $(M^n,\cF)$ be a foliation of dimension $k$. On the manifold
$M^n\times \bR$ we define naturally a foliation $\wh{\cF}$ of
dimension $k+1$ whose leaves are manifolds $L_\al\times {\bR}$,
$\al\in A$, where $L_\al$, $\al\in A$, are the leaves of the
foliation $\cF$.

%        LEMMA 2.3.
\begin{lem}
The projection $\pi:M^n\times {\bR}\to M^n$ defines an isomorphism
in vertical cohomologies.
\end{lem}

%        PROOF.
\begin{pf}
Consider the zero section $s$ of the trivial bundle $M^n\times
{\bR} @>{\pi}>> M^n$, i.e., $s(x)=(x,0)$, $x\in M^n$. Then the
mappings $\pi$ and $s$ are the leaf-to-leaf transforming mappings
which define the cochain mappings
$\pi^*:(\Lb_{\cF}^*(M^n),d_{\cF}^*) \to
    (\Lb_{\wh{\cF}}^*(M^n\times {\bR}),d_{\wh{\cF}}^*)$ and
$s^*:(\Lb_{\wh{\cF}}^*(M^n\times {\bR}),d_{\wh{\cF}}^*) \to
    (\Lb_{\cF}^*(M^n),d_{\cF}^*)$. Since $s^*\circ \pi^*=1$,
$\pi^*$ is an monomorphism. We shall show that $\pi^*$ induces an
isomorphism at the cohomology level. To this end, we shall
construct a cochain equivalence of the mappings 1 and $\pi^*\circ
s^*$.

        Note that each form from $\Lb_{\cF}^*(M^n\times {\bR})$
can be uniquely represented by linear combinations of the
following two types of forms:
\begin{align*}
        \text{(I)} \;\; & (\pi^*\vf)\cdot f, \;\;\; \vf\in \Lb_{\cF}^*(M^n),
        \;\;\;f\in C^\infty(M^n\times {\bR});  \\
        \text{(II)} \;\; & (\pi^*\vf)\wedge dt\cdot f, \;\;\;
        \vf\in \Lb_{\cF}^*(M^n),
        \;\;\;f\in C^\infty(M^n\times {\bR}),
\end{align*}
where $t$ is the coordinate on the straight line ${\bR}$. Define
the operator $K^*:\Lb_{\wh{\cF}}^*(M^n\times {\bR})\to
{\Lb_{\wh{\cF}}^*}^{-1}(M^n\times {\bR})$ as follows:
\begin{align*}
    & K^*((\pi^*\vf)\cdot f)=0,  \\
    & K^*((\pi^*\vf)\wedge dt\cdot f)=\pi^*\vf\cdot \int_0^t f\,dt.
\end{align*}

        A direct calculation shows that the relation
$$  1-\pi^*\circ s^*=(-1)^{q-1} (d_{\wh{\cF}}^{q-1}K^q-K^{q+1}d_{\wh{\cF}}^q)   $$
is fulfilled on the forms of types (I) and (II). %Thus the lemma is complete.
\end{pf}

%        DEFINITION 2.4.
\begin{defi}
Two leaf-to-leaf transforming mappings $f,g:({M_1}^n,{\cF}_1)
\linebreak  \to ({M_2}^m,{\cF}_2)$ between the foliations
$({M_1}^n,{\cF}_1)$ and $({M_2}^m,{\cF}_2)$ are called
leaf-to-leaf transforming homotopic if there exists a leaf-to-leaf
transforming mapping
$$  F:(M_1^n\times {\bR},\wh{\cF}_1)\to (M_2^m,{\cF}_2)     $$
such that
$$  \begin{cases}
        F(x,t)=f(x), & t\geq 1, \\
        F(x,t)=g(x), & t\leq 0,
    \end{cases} \quad x\in M,\;\;\;t\in {\bR}.        $$
\end{defi}

%        THEOREM 2.5
\begin{th}[A Homotopy Axiom]
Leaf-to-leaf transforming  ho- \linebreak motopic mappings induce
identical map\-pings in vertical cohomologies.
\end{th}

%        PROOF.
\begin{pf}
Let $f,g:(M_1^n,{\cF}_1)\to ({M_2}^m,{\cF}_2)$ be the leaf-to-leaf
transforming homotopic mappings and $F:(M_1^n\times
{\bR},\wh{\cF}_1)\to (M_2^m,{\cF}_2)$ be the homotopy between $f$
and $g$. Denote by $s_0$ and $s_1$ the sections
$s_0,s_1:(M_1^n,{\cF}_1)\to (M_1^n\times {\bR},\wh{\cF}_1)$,
$s_0(x)=(x,0)$, $s_1(x)=(x,1)$, $x\in M_1^n$. Then $f=F\circ s_1$
and $g=F\circ s_0$. Hence we have $f^*=s_1^*\circ F^*$ and
$g^*=s_0^*\circ F^*$. From the proof of Lemma 2.3 it follows that
$s_1^*=s_0^*=(\pi_1^*)^{-1}$, where $\pi_1:M_1^n\times {\bR}\to
M_1^n$
is the projection. Therefore $f^*=g^*$. %Thus the theorem is complete.
\end{pf}

        The foliations $(M_1^n,{\cF}_1)$ and
$(M_2^m,{\cF}_2)$ will be said to be  of the same homotopy type if
there are  leaf-to-leaf transforming smooth mappings $f:M_1^n\to
M_2^m$ and $g:M_2^m\to M_1^n$ such that $g\circ f$ and $f\circ g$
are leaf-to-leaf transforming homotopic to the identical mappings
of the foliations $(M_1^n,{\cF}_1)$ and $(M_2^m,{\cF}_2)$,
respectively.

%        COROLLARY 2.6.
\begin{cor}
If two foliations $(M_1^n,{\cF}_1)$ and $(M_2^m,{\cF}_2)$ are of
the same homotopy type, then their vertical cohomologies are
isomorphic.
\end{cor}

%        2.3.
\subsection{Relative Vertical Cohomologies}

        Let $(M^n,{\cF}_1)$ and $(N^m,{\cF}_2)$ be two foliations, and let
$f$ be a leaf-to-leaf transforming smooth mapping $f:M^n\to N^m$.
Define the differential complex
$$  (\Lb^*(f),\ol{d}{}^*), \quad \Lb^*(f)=
        \us{q\geq 0}{\oplus} \Lb^q(f),      $$
where
$$  \Lb^q(f)=\Lb_{{\cF}_2}^q(N^m)\oplus \Lb_{{\cF}_1}^{q-1}(M^n), \quad
        \ol{d}{}^*(\om,\tht)=(-d_{{\cF}_2}^*\om,
            f^*\om+d_{{\cF}_1}^*\tht).    $$

        We easily verify that $\ol{d}{}^2=0$ and denote the cohomology
groups of this complex by $H^*(f)$. Note that the complex
$(\Lb^*(f),\ol{d}{}^*)$ is the cone of the cochain mapping
$f^*:\Lb_{{\cF}_2}^*(N^m)\to \Lb_{{\cF}_1}^*(M^n)$. If we
regraduate the complex $\Lb_{{\cF}_1}^*(M^n)$ as
$\wt{\Lb}_{{\cF}_1}^p(M_1^n)\equiv \Lb_{{\cF}_1}^{p-1}(M^n)$, then
we obtain an exact sequence of differential complexes
%   (2.5)
\begin{equation}
    0 @>>> \wt{\Lb}_{{\cF}_1}^*(M^n) @>{\al}>> \Lb^*(f) @>{\bt}>>
        \Lb_{{\cF}_2}^*(N^m) @>>> 0
\end{equation}
with the obvious mappings $\al$ and $\bt$: $\al(\tht)=(0,\tht)$,
$\bt(\om,\tht)=\om$. From (2.5) we have an exact sequence in
cohomologies
$$  \cdots @>>> H_{{\cF}_1}^{q-1}(M^n) @>{\al^*}>> H^q(f) @>{\bt^*}>>
        H_{{\cF}_2}^q(N^m) @>{\dl^*}>> H_{{\cF}_1}^q(M^n)
        @>>> \cdots\;.      $$

        It is easily seen that $\dl^*=f^*$. Let $\om\in \Lb_{{\cF}_2}^q(N^m)$
be the closed form, and $(\om,\tht)\in \Lb^q(f)$. Then
$\ol{d}(\om,\tht)=(0,f^*\om+d_{{\cF}_1}^{q-1}\tht)$, and by the
definition of the operator $\dl^*$ we have
$\dl^*[\om]=[f^*\om+d_{{\cF}_1}^{q-1}\tht]=f^*[\om]$. Hence we
finally get a long exact sequence
%   (2.6)
\begin{equation}
\ \hskip-1cm    \cdots @>>> H_{{\cF}_1}^{q-1}(M^n) @>{\al^*}>>
H^q(f) @>{\bt^*}>>
        H_{{\cF}_2}^q(N^m) @>{f^*}>> H_{{\cF}_1}^q(M^n)
        @>{\al^*}>> \cdots \;.
\end{equation}

%        COROLLARY 2.7.
\begin{cor}
If the foliations $(M^n,{\cF}_1)$ and $(N^m,{\cF}_2)$ are of the
$p$-th and $q$-th dimension, respectively, then
\begin{itemize}
\item[(i)] $\bt^*:H^{p+1}(f)\to H_{{\cF}_2}^{p+1}(N^m)$
        is an epimorphism, \\
            $\al^*:H_{{\cF}_1}^q(M^n)\to H^{q+1}(f)$ is an epimorphism, \\
            $\bt^*:H^i(f)\to H_{{\cF}_2}^i(N^m)$ is an isomorphism for $i>p+1$, \\
            $\al^*:H_{{\cF}_1}^i(M^n)\to H^{i+1}(f)$ is an isomorphism for $i>q$;

\item[(ii)] $H^i(f)=0$ for $i>\max\{p+1,q\}$.
\end{itemize}
\end{cor}

 %       THEOREM 2.8.
\begin{th}
If $f,g:(M^n,{\cF}_1)\to (N^m,{\cF}_2)$ are leaf-to-leaf
transforming homotopic mappings, then $H^*(f)=H^*(g)$.
\end{th}

%        PROOF.
\begin{pf}
Let $F:(M^n\times {\bR},\wh{{\cF}}_1)\to (N^m,{\cF}_2)$ be the
homotopy mapping between $f$ and $g$. Let $s_0$ and $s_1$ be the
zero and the unit section, respectively, of the trivial bundle
$M^n\times {\bR}@>{\pi}>> M^n$. Then $F\circ s_0=g$ and $F\circ
s_1=f$. Hence we have a homomorphism between the short exact
sequences
$$  \begin{CD}
        0 @>>> \wt{\Lb}_{\wh{{\cF}}_1}^*(M^n\times {\bR})
            @>{\al'}>> \Lb^*(F) @>{\bt'}>>
            \Lb_{{\cF}_2}^*(N^m) @>>> 0 \\
        @. @VV{s_1^*}V @VV{id\oplus s_1^*}V @VV{id}V \\
        0 @>>> \wt{\Lb}_{{\cF}_1}^*(M^n) @>>> \Lb^*(f) @>>>
            \Lb_{{\cF}_2}^*(N^m) @>>> 0
    \end{CD}\;.     $$

        This homomorphism defines a homomorphism  between the corresponding
long cohomologic sequences
$$  \begin{smallmatrix}%array}{ccccccccccccc}
        \cdots & \to & H_{{\cF}_2}^q(N^m) & \to
            & H_{\wh{{\cF}}_1}^q(M^n\times {\bR}) & \to
            & H^{q+1}(F) & \to
            & H_{{\cF}_2}^{q+1}(N^m) & \to
            & H_{\wh{{\cF}}_1}^{q+1}(M^n\times {\bR}) & \to
            & \cdots \\
        & & \downarrow{\vcenter{\rlap{$\scriptstyle id$}}}
            & & \downarrow{\vcenter{\rlap{$\scriptstyle s_1^*$}}}
            & & \downarrow{\vcenter{\rlap{$\scriptstyle \gm$}}}
            & & \downarrow{\vcenter{\rlap{$\scriptstyle id$}}}
            & & \downarrow{\vcenter{\rlap{$\scriptstyle s_1^*$}}}&&\\
        \cdots & \to & H_{{\cF}_2}^q(N^m) & \to
            & H_{{\cF}_1}^q(M^n) & \to
            & H^{q+1}(f) & \to
            & H_{{\cF}_2}^{q+1}(N^m) & \to
            & H_{{\cF}_1}^{q+1}(M^n) & \to
            & \cdots
    \end{smallmatrix} \;, $$%array}\;,          $$
where $\gm$ is the mapping induced by $id\oplus s_1^*$. Since
$s_1^*$ is an isomorphism (Lemma 2.3), by virtue of the lemma on
five homomorphisms we conclude that $\gm$ is also an isomorphism,
i.e., $H^*(f)\approx H^*(F)$. By a similar reasoning we can
conclude that $H^*(g)\approx H^*(F)$.
\end{pf}

        If $(M^n,{\cF}_1)$ is a subfoliation of  the  foliation
$(W^m,{\cF}_2)$, i. e.,  the embedding $M^n \str{j}{\hook} W^m$ is
simultaneously a leaf-to-leaf transforming mapping, then the
cohomology algebra $H^*(j)$ will be said to be the algebra of
relative vertical cohomologies. Denote it by
$H_{{\cF}_2,{\cF}_1}^*(W;M)$.

        Now sequence (2.6) can be rewritten as
%   (2.7)
\begin{gather}
    \cdots @>>> H_{{\cF}_1}^{q-1}(M) @>{\al^*}>> H_{{\cF}_2,{\cF}_1}^q(W;M)
        @>{\bt^*}>> H_{{\cF}_2}^q(W) @>{j^*}>> \notag  \\
    @>{j^*}>> H_{{\cF}_1}^q(M)
        @>{\al^*}>> \cdots \;.
\end{gather}

        Note that if we forget the structure of the foliation,
then, as is known, the embedding $M^n \str{j}{\hook} W^m$ defines
a long exact cohomological sequence of the pair $(W^m,M^n)$ in de
Rham's theory. One can easily verify that the homomorphism $i^*$
from Proposition 2.2 defines a morphism between the long exact
cohomological sequence of the pair $(W^m,M^n)$ in de Rham's theory
and sequence (2.7).

%        2.4.
\subsection{The Generalized Mayer--Vietoris Principle for Vertical Cohomologies.
            A Combinatorial Definition of Vertical Cohomologies.}
        Let  $(M^n,{\cF})$ be a smooth foliation of dimension
$k$; let $U=\{u_\al\}_{\al\in A}$ be an open countable covering of
the manifold $M^n$. Similarly to \v{C}ech-de Rham's theory, we
define a double complex which will be used to calculate vertical
cohomologies of the foliation $(M^n,{\cF})$.

        Denote by $u_{\al_0\cdots\al_p}$ the intersection of open sets
$u_0,\dots,u_p$ and by $\coprod$ the disjunctive union. Then we
have a sequence of open sets
$$  M^n \;\longleftarrow \;\coprod_{\al_0} u_{\al_0}\;
            \begin{matrix}
                {\scriptstyle \pa_0} \\[-2.5mm]
                \longleftarrow \\[-3mm]
                \longleftarrow \\[-2.5mm]
                {\scriptstyle \pa_1}
            \end{matrix} \;
        \coprod_{\al_0<\al_1} u_{\al_0\al_1} \;
            \begin{matrix}
                {\scriptstyle \pa_0} \\[-2.5mm]
                \longleftarrow \\[-3mm]
                {\scriptstyle \pa_1} \\[-3mm]
                \longleftarrow \\[-3mm]
                \longleftarrow \\[-2.5mm]
                {\scriptstyle \pa_2}
            \end{matrix} \;
        \coprod_{\al_0<\al_1<\al_2} u_{\al_0\al_1\al_2} \;
            \begin{matrix}
                \longleftarrow \\[-3mm]
                \longleftarrow \\[-3mm]
                \longleftarrow \\[-3mm]
                \longleftarrow
            \end{matrix} \; \cdots\;,           $$
where $\pa_i$ is the embedding $\pa_i(u_{\al_0\cdots
\al_p})=u_{\al_0\cdots \wh{\al}_i\cdots\al_p}$. This sequence of
embeddings induces a sequence of restriction mappings of vertical
forms
\begin{gather*}
    \Lb_{{\cF}}^*(M^n) \; \str{r^*}{\longrightarrow} \;
        \prod_{\al_0} \Lb_{{\cF}}^*(u_{\al_0})\;
            \begin{matrix}
                {\scriptstyle \dl_0} \\[-2.5mm]
                \longrightarrow \\[-3mm]
                \longrightarrow \\[-2.5mm]
                {\scriptstyle \dl_1}
            \end{matrix} \;
        \prod_{\al_0<\al_1} \Lb_{{\cF}}^*(u_{\al_0\al_1}) \;
            \begin{matrix}
                {\scriptstyle \dl_0} \\[-2.5mm]
                \longrightarrow \\[-3mm]
                {\scriptstyle \dl_1} \\[-3mm]
                \longrightarrow \\[-3mm]
                \longrightarrow \\[-2.5mm]
                {\scriptstyle \dl_2}
            \end{matrix} \;  \\
        \begin{matrix}
            {\scriptstyle \dl_0} \\[-2.5mm]
            \longrightarrow \\[-3mm]
            {\scriptstyle \dl_1} \\[-3mm]
            \longrightarrow \\[-3mm]
            \longrightarrow \\[-2.5mm]
            {\scriptstyle \dl_2}
        \end{matrix} \;
        \prod_{\al_0<\al_1<\al_2} \Lb_{{\cF}}^*(u_{\al_0\al_1\al_2}) \;
            \begin{matrix}
                \longrightarrow \\[-3mm]
                \longrightarrow \\[-3mm]
                \longrightarrow \\[-3mm]
                \longrightarrow
            \end{matrix} \; \cdots\;,
\end{gather*}
where $\dl_i$ is induced by the imbedding $\pa_i$, i.e.,
$\dl_i=(\pa_i)^*$.

        Let us define the difference operator $\dl$ by the following
rule: if $\om_{\al_0\cdots\al_p}\in
\Lb_{{\cF}}^q(u_{\al_0\cdots\al_p})$ denotes the components of the
element $\om\!\in\!\! \prod\limits_{\al_0<\cdots<\al_p}\!\!\!\!
\Lb_{{\cF}}^q(u_{\al_0\cdots\al_p})$, then
%   (2.8)
\begin{equation}
    (\dl\om)_{\al_0\cdots\al_{p+1}}=\sum_{i=0}^{p+1}(-1)^i
        \om_{\al_0\cdots\wh{\al}_i\cdots\al_{p+1}},
\end{equation}
where $\om_{\al_0\cdots\wh{\al}_i\cdots\al_{p+1}}\equiv \dl_i
        (\om_{\al_0\cdots\al_{i-1}\al_{i+1}\cdots\al_{p+1}})$.

        By a standard reasoning one can verify that $\dl^2=0$.

        Consider now the double complex
%   (2.9)
\begin{equation}
    K^{p,q}\equiv C^p(U,\Lb_{{\cF}}^q)=\prod_{\al_0<\cdots<\al_p}
        \Lb_{{\cF}}^q(u_{\al_0\cdots\al_p}).
\end{equation}

whose horizontal mappings are the operators $\dl^*$ and vertical
mappings are the operators $d_{{\cF}}^*$. As is known, this double
complex can be reduced to an ordinary differential complex
$(K^*,D^*)$:
$$  K^n=\us{p+q=n}{\oplus} K^{p,q}, \quad D^n=\dl^p+(-1)^pd_{{\cF}}^q
        \quad \text{on} \quad K^{p,q}.      $$

%        LEMMA 2.9.
\begin{lem}
The sequence
$$  \cdots @>>> \Lb_{{\cF}}^*(M^n) @>{r^*}>> K^{0,*} @>{\dl^0}>> K^{1,*}
        @>{\dl^1}>> K^{2,*} @>{\dl^2}>> \cdots      $$
is exact.
\end{lem}

%        PROOF.
\begin{pf}
Let $q\geq 0$ be an integer number. Clearly, $\Lb_{{\cF}}^q(M^n)$
is the kernel of the first operator $\dl^0$, since an element from
$\prod\limits_{\al_0}\Lb^q(u_{\al_0})$ is a global form on $M^n$
if and only if its components consistent at the intersections.

        Let $\{\tht_\al\}_{\al\in A}$ be the partitioning of unity
subordinate to the covering $U=\{u_\al\}_{\al\in A}$. If $\om\in
K^{p,q}$ is the cocycle of the operator $\dl^p$, then we can
assign to it $p-1$ cochains $K\om$ by the formula
$(K\om)_{\al_0\cdots\al_{p-1}}=\sum\limits_\al \tht_\al
    \om_{\al\al_0\cdots\al_{p-1}}$ (it is assumed here that
$\om_{\dots,\al,\dots,\bt,\dots}=-\om_{\dots,\bt,\dots,\al,\dots}$,
if $\al>\bt$; clearly, this is consistent with the operation
$\dl$, i.e., $(\dl\om)_{\dots,\bt,\dots,\al,\dots}=
    -(\dl\om)_{\dots,\al,\dots,\bt,\dots}$). In that case
\begin{gather*}
    (\dl^{p-1}K\om)_{\al_0\cdots\al_p}=\sum_i (-1)^i
        (K\om)_{\al_0\cdots\wh{\al}_i\cdots\al_p}=
        \sum_{i,\al} (-1)^i\tht_\al
        \om_{\al\al_0\cdots\wh{\al}_i\cdots\al_p}= \\
    =\sum_\al\tht_\al[\om_{\al_0\cdots,\al_p}-
        (\dl\om)_{\al\al_0\cdots\al_p}]=
        \sum_\al\tht_\al \om_{\al_0\cdots\al_p}=\om_{\al_0\cdots\al_p}.
\end{gather*}

        Hence $\dl^{p-1}(K\om)=\om$.
\end{pf}

%        THEOREM 2.10.
\begin{th}
The cohomologies of the double complex $K^{p,q}$, i.e., the
cohomologies of the complex $(K^*,D^*)$ are isomorphic to the
vertical cohomologies $H_{{\cF}}^*(M^n)$. This isomorphism is
obtained  by the mapping of the restriction $r^*$.
\end{th}

%        PROOF.
\begin{pf}
Since
$D^*r^*=(\dl^0+d_{\cF}^*)r^*=\dl^0r^*+d_{\cF}^*r^*=d_{\cF}^*r^*=
    r^{*+1} d_{\cF}^*$, $r^*$ is a cochain mapping, it induces the
mapping in cohomologies $r_*$.

        Let $\vf\in K^m$, be the cocycle, i.e., $D^m\vf=0$. We can represent
$\vf$ as a sum $\vf=\vf^{0,m}+\vf^{1,m-1}+\cdots+\vf^{p,m-p}$,
where $\vf^{i,j}\in K^{ij}$ and $\vf^{p,m-p}\neq 0$. Then
$\dl^p\vf^{p,m-p}=0$. Hence by Lemma 2.9 we find that there exists
an element $\vf'\in K^{p-1,m-p}$ such that
$\dl^{p-1}\vf'=\vf^{p,m-p}$. Then the element $\vf-D^{m-1}\vf'$ is
obviously cohomologic to $\vf$ and has no component in
$K^{p,m-p}$. After repeating this procedure several times, we
obtain an element $\wt{\vf}\in K^{0,m}$ which is cohomologic to
$\vf$. For this element we have $d_{\cF}^m\wt{\vf}=0$,
$\dl^0\wt{\vf}=0$, i.e., $\wt{\vf}$ defines the global vertical
closed form on $M^n$ which by means of $r^m$ transforms to
$\wt{\vf}$. This shows that the mapping $r_*^m$ is epimorphic. Let
us now show that the mapping $r_*^m$ is monomorphic.

        Let $\vf\in \Lb_{\cF}^m(M^n)$ and $d_{\cF}^m\vf=0$. Then if
$r^m(\vf)=D^{m-1}\vf'$, $\vf'\in K^{m-1}$,
$\vf'=\vf'{}^{0,m-1}+\cdots+\vf'{}^{p,m-p-1}$, where
$\vf'{}^{p,m-p-1}\neq 0$, then $\dl^p\vf'{}^{p,m-p-1}=0$ (because
$D^{m-1}\vf'$ has only one component in $K^{0,m}$).

        Proceeding as above, we can find an element $\vf''$ such that
$\vf''\in K^{0,m-1}$ and $\vf'-\vf''\in D^{m-2}(K^{m-2})$. Then
$r^m(\vf)=d_{\cF}^{m-1}\vf''$ and $\dl^0\vf''=0$. Therefore
$\vf''$ defines a global form $\wt{\vf}''$ on $M^n$ such that
$d_{\cF}^{m-1}\wt{\vf}''=\vf$.
\end{pf}  %        Thus the proof is complete.

        By analogy with de Rham's theory Theorem 2.10 can be called
the generalized Mayer--Vietoris principle.

        From the lower row of the double complex $K^{*,*}$ let us
choose a subcomplex, namely, a kernel of the differential
$d_{\cF}^0$. We get a sequence
%   (2.10)
\begin{equation}
    C_{\cF}^0(U) @>{\dl^0}>> C_{\cF}^1(U) @>{\dl^1}>>
        C_{\cF}^2(U) @>{\dl^2}>> \cdots \;,
\end{equation}
where $C_{\cF}^i(U)\equiv \ker d_{\cF}^0\sbs
\prod\limits_{\al_0<\cdots<\al_i}
    \Lb_{\cF}^0(u_{\al_0\cdots \al_i})$, and $\dl^*$
is the difference operator defined above. The cohomologies of
complex (2.10) will be called the \v{C}ech cohomologies of the
foliation $(M^n,\cF)$ for the covering $U$. They are a purely
combinatorial object and will be denoted by $H_{\cF}^*(U)$.  Let
the covering $U=\{u_\al\}_{\al\in A}$ consist of foliated open
sets such that all finite nonempty intersections are contractible.
Any manifold is known to have such a covering. Then, in view of
the fact that Poincare's lemma ([I]) is valid for vertical
cohomologies, we obtain

%        LEMMA 2.11.
\begin{lem}
The sequence
$$  0 @>>> C_{\cF}^p(U) @>{j^p}>> \prod_{\al_0<\cdots<\al_p}
        \Lb_{\cF}^0 (u_{\al_0\cdots\al_p})
        @>{d_{\cF}^0}>> \prod_{\al_0<\cdots<\al_p}
        \Lb^1 (u_{\al_0\cdots\al_p}) @>{d_{\cF}^1}>> \cdots     $$
is exact, $p\geq 0$, where $j^p$ is the embedding.
\end{lem}

By Lemma 2.11 and the same reasoning as we used in proving Theorem
2.10  we can prove

%        THEOREM 2.12.
\begin{th}
The cochain mapping $j^*$ defines an isomorphism between the
\v{C}ech cohomologies $H_{\cF}^*(U)$ and cohomologies of the
double complex $K^{*,*}$. Hence the \v{C}ech cohomologies and
vertical ones are isomorphic:
$$  H_{\cF}^*(U)\approx H_{\cF}^*(M^n). $$
\end{th}

%        COROLLARY 2.13.
\begin{cor}
A zero-dimensional vertical cohomology $H_{\cF}^0(M^n)$ is
\linebreak isomorphic to a group of smooth functions on the
foliated manifold  $M^n$, which are constant on the leaves.
\end{cor}

%        \S 3.
\section{Completely Integrable Hamiltonian Systems and Vertical
        Cohomologies}

%        3.1.
\subsection{Topology of Constant Energy Surfaces of Completely Integrable
                Hamiltonian Systems}

        In this subsection we shall briefly recall the basic facts from the
topological theory of integrable Hamiltonian systems ([5]).

        Let $v=\sgrad H$ be a Hamiltonian system on the symplectic manifold
$(M^{2n},\om)$, where $\om$ is a symplectic 2-form on $M^n$, and
let $v$ be integrable. Thus there exist $n$ independent (almost
everywhere) smooth integrals $f_1=H$, $f_2,\dots,f_n$ in
involution, i.e., $\{f_i,f_j\}=0$, $i,j=1,n$, where $\{\ ,\ \}$ is
the Poisson bracket. Let $F:M^{2n}\to {\bR}^n$ be the moment
mapping which corresponds to these integrals, i.e.,
$F(x)=(f_1(x),\dots,f_n(x))$, $x\in M^{2n}$. Let $N$ be the set of
critical points of the moment mapping,
 and $\Sigma=F(N)$ be the set of all critical values which is
called the bifurcation diagram.

        Clearly, we have two cases: (a) $\dim \Sigma<n-1$ and
(b)~$\dim\Sigma=n-1$. In the case (a) the set $\Sigma$ does not
separate the space ${\bR}^n$ and therefore  all nonsingular leaves
$B_a=F^{-1}(a)$ are diffeomorphic to one another (it is well known
that if they are compact, then they are diffeomorphic to the tori
$T^n$, and if they are noncompact, then they are diffeomorphic to
the cylinders $T^k\times {\bR}^{n-k}$). The case (b) is more
difficult. Below we shall consider a theorem from [5].

        Suppose that the restriction $f=f_1|_{X^{n+1}}$ to a joint compact
nonsingular surface of the level of the rest of $n-1$ integrals
$X^{n+1}=\{x\in M^{2n}|f_i(x)=c_i$, $i=\ol{2,n}\}$ is a Bott
function, i.e., all critical points of this restriction are
organized into nondegenerate critical submanifolds (a critical
submanifold $L^k\sbs X^{n+1}$ is nondegenerate if the restriction
of the function $f$ to every normal plane $P^{n+1-k}$ has a
nondegenerate Morse singularity at the point $P^{n+1-k}\cap L^k$).

        Let $c_1$ be a critical value of the function $f$ on
the surface $X^{n+1}$, and let $c=(c_1,\dots,c_n)$, i.e.,
$c=(c_1,c_2,\dots,c_n)\in \Sigma$. Let $B_c=F^{-1}(c)$ be a
critical fiber of the moment mapping. Thus $B_c=\{f_1=c_1\}$ is a
critical level surface of $f_1$ on $X^{n+1}$.

%        THEOREM 3.1
\begin{th}[A. T. Fomenko]
Each connected compact component $B_c^0$ of the critical fiber
$B_c$ is homeomorphic to a set which is one of the following four
types: $(1)$ a torus $T^n$; $(2)$ the nonorientable manifolds
$K_0^n$ and $K_1^n$; $(3)$ a torus $T^{n-1}$; or $(4)$ a cell
complex $T^n_1\cup T^n_2$ obtained by removing  $n-1$-dimensional
tori $T^{n-1}_1$ from $T_1^n$ and $T^{n-1}_2$ from $T_2^n$, which
realize nonzero generators of the homology groups
$H_{n-1}(T_1^n,{\bZ})$ and $H_{n-1}(T_2^n,{\bZ})$,  and glueing
$T_1^n$ and $T_2^n$ together by identifying  only the tori
$T_1^{n-1}$ and $T_2^{n-1}$ $($by means of a diffeomorphism$)$. In
cases $(1)$--$(3)$  the critical fibers  consist entirely of
critical points of the function $f$ on which a maximum or a
minimum is attained. In case $(4)$ the critical points of $f$ in
the critical fiber $B^0_c$ forms a torus $T^{n-1}$ $($the result
of glueing together the tori $T^{n-1}_1$ and $T^{n-1}_2)$ which is
a ``saddle'' for the function $f$.
\end{th}

        The manifolds $K_1^n$ and $K_2^n$ from Theorem 3.1 are the factor-sets
of the torus $T^n$ generated by the action of some group on $T^n$
without fixed points ([5]).

    It is important to investigate a particular case of
the integrable Hamiltonian system $v=\sgrad H$ on the
four-dimensional symplectic manifold $(M^4,\om)$. For mechanical
and physical reasons, it is useful to study the integrability
effect on an individual isoenergetic surface $Q^3\sbs M^4$ given
by the equation $H(x)=h$, where $h$ is a regular value of the
Hamiltonian $H$. The restriction of the system $v$ on $Q^3$ will
be denoted by the same symbol $v$. In what follows $Q^3$ will be
assumed to be the closed manifold. We also assume that an
additional integral $f$ is the Bott function on $Q^3$. Then to the
critical fibers from Theorem 3.1 there correspond the following
manifolds:

(1)~the torus $T^2$;

(2)~Klein's bottle $K^2$;

(3)~the circle $S^1$; (4) a piecewise smooth two-dimensional
polyhedron with a singularity of the type of "a fourfold line" (a
transversal intersection of two planes).

        The Hamiltonian $H$ is called the nonresonance one on $Q^3$
if the set of Liouville tori with irrational windings is
everywhere dense in $Q^3$.

    We say that the Hamiltonian system $v$ is integrable in the Bott sense
if among Bott functions there is an additional first integral of
$v$.

%        DEFINITION 3.2
\begin{defi}[{[6]}]
Two  integrable, in the Bott sense, nonresonance Hamiltonian
systems $v_1,v_2$ on the oriented manifolds $Q_1^3,Q_2^3$ are said
to be topologically equivalent if there exists an orientation
preserving diffeomorphism $g:Q_1^3\to Q_2^3$ which transforms the
Liouville tori of the system $v_1$ to those of the system $v_2$
(critical tori are included into the number of Liouville tori),
and the isolated critical circles to the isolated critical circles
with the same orientation which is by the field $v$.
\end{defi}

The partitioning of the manifold $Q^3$ into Liouville tori and
critical level surfaces of the Bott integral $f$ is called the
Liouville foliation on $Q^3$ (for a given integrable nonresonance
Hamiltonian system $v$). Obviously, the diffeomorphism $g$,
appearing in Definition 3.2, is a leaf-to-leaf mapping between
Liouville foliated manifolds.

        Note that the Liouville foliation defined above does not depend on
a choice of the Bott integral $f$.

        An analogous definition can be introduced for the multidimensional
case as well.  Two completely integrable, in a Bott sense,
nonresonance Hamiltonian systems $v_1$ and $v_2$ on the
nonsingular level surfaces $X_1^{n+1}$ and $X_2^{n+1}$ are said to
be topologically equivalent if there exists a diffeomorphism
$h:X_1^{n+1}\to X_2^{n+1}$ which transforms the Liouville tori of
the system $v_1$ to those of the system $v_2$ (critical tori $T^n$
are included into the number of Liouville tori), and critical
fibers $B_c$ of the system $v_1$ to those of the system $v_2$.

Note that the Liouville foliation is not a foliation in the sense
of \S 2.

Let $v=\sgrad H$ be a completely integrable, in the Bott sense,
nonresonance Hamiltonian system on the isoenergetic surface $Q^3$.
As we already know, this system defines the Liouville foliation on
$Q^3$. Denote it by $\cF$.

        If from $Q^3$ we remove all critical
(i. e., maximal, minimal and "saddle"-type)  circles $S^1$, then
on the resulting open manifold $Q^{'3}$  the foliation $\cF$
induces a true foliation $\cF'$. Thus we have defined the vertical
cohomology groups $H_{\cF'}^i(Q^{'3})$, $i=0,1,2$. Similarly, one
can derive the vertical cohomology groups $H_{\cF'}^i(X^{'n+1})$
in the multidimensional case, $i=0,1,\dots,m$.

        If now  $v_1$ and $v_2$ are assumed to be the topologically
equivalent systems on the manifolds $Q_1^3$ and $Q_2^3$,
respectively, then the diffeomorphism $g:Q_1^3\to Q_2^3$ which
preserves the Liouville foliation obviously induces the
leaf-to-leaf transforming diffeomorphism $g':Q_1^{'3}\to
Q_2^{'3}$. Therefore the following theorem is valid.

%        THEOREM 3.3.
\begin{th}
The cohomology groups  $H_{\cF'}^i(Q^{'3})$ are a topological
invariant of completely integrable, in the Bott sense,
nonresonance Hamiltonian systems $v=\sgrad H$ on the isoenergetic
surface $Q^3$.
\end{th}

       Note that in the multi-dimensional case the groups
$H_{\cF'}^i(X^{'n+1})$, $i=\ol{0,n}$, are also  topological
invariants.

      Thus to determine a topological invariant we had to deal with
the open manifold $Q^{'3}$. To avoid this, we shall slightly
modify the definition of cohomology groups.

       The isoenergetic surface $Q^3$ can be considered  as the set of
Liouville tori, circles $S^1$ (which are maximal, minimal and
"saddle"-type)  and rings of the type $S^1\times (0,1)$ (which are
obtained from a hyperbolic critical fiber of type (4) (Theorem
3.1) by removing hyberbolic critical circles of the integral $f$).
We denote this partitioning of the manifold $Q^3$ by $\cP$. By
$S(Q^3)$ we denote a Lie algebra of smooth vector fields on $Q^3$
and by $\ol{S}(Q^3)$ a Lie subalgebra of the algebra $S(Q^3)$
consisting of smooth vector fields on $Q^3$ tangent to the
manifolds of the partitioning $\cP$. The algebra $C^\infty(Q^3)$
is, obviously, a module over $\ol{S}(Q^3)$ with respect to the
product $x\cdot g=xg$, $x\in \ol{S}(Q^3)$, $g\in C^\infty(Q^3)$.

        Under a $k$-dimensional cochain of the algebra $\ol{S}(Q^3)$
with coefficients in $C^\infty(Q^3)$ we shall mean an element of
the space $A^k(\ol{S}(Q^3))$ which is a skew-symmetric $k$-linear
functional on $\ol{S}(Q^3)$ with values in $C^\infty(Q^3)$. The
differential $d^k:A^k(\ol{S}(Q^3))\to A^{k+1}(\ol{S}(Q^3))$ is
defined by the formula
\begin{gather*}
    d^k\al(X_1,\dots,X_{k+1})=\sum_{i=1}^{k+1}(-1)^{i-1}X_i
        \al(X_1,\dots,\wh{X}_i,X_{i+1},\dots,X_{k+1})+  \\
    +\sum_{i<j} (-1)^{i+j} \al([X_i,X_j],X_1,\dots,\wh{X}_i,\dots,\wh{X}_j,
        \dots,X_{k+1}),
\end{gather*}
where $\al\in A^k(\ol{S}(Q^3))$; $X_1,\dots,X_{k+1}\in
\ol{S}(Q^3)$.

        It is easy to verify that $d^k\circ d^{k-1}=0$. Thus
$(A^*(\ol{S}(Q^3)),d^*)$ is a cochain complex. We denote the
cohomologies of the complex $(A^*(\ol{S}(Q^3)),d^*)$ by
$H_{\cP}^*(Q^3)$. Obviously, $H_{\cP}^i(Q^3)=0$, $i>2$.

        We have also defined the cohomology groups $H_{\cP}^*(X^{n+1})$
for the miltidimensional case. Note that here
$H_{\cP}^i(X^{n+1})=0$, $i>n$.

        Obviously, $A^*(\ol{S}(Q^3))$ is a graduated ring, and
the differential $d$ is the antidifferentiation. Hence the
cohomologies $H_{\cP}^*(Q^3)$ also acquire a ring structure.

        If now $v_1$ and $v_2$ are as above the topologically equivalent
systems on the manifolds $Q_1^3$ and $Q_2^3$, respectively, then
the diffeomorphism $g:Q_1^3\to Q_2^3$ defines an isomorphism
between the differential complexes $(A^*(\ol{S}(Q_1^3)),d^*)$ and
$(A^k(\ol{S}(Q_2^3)),d^*)$. Therefore we have

%        THEOREM 3.4.
\begin{th}
Let $v=\sgrad H$ be an integrable, in the Bott sense, nonresonance
Hamiltonian system on the isoenergetic surface $Q^3$, and let
$\cP$ be the partitioning of the manifold $Q^3$ into Liouville
tori, critical circles of the additional Bott integral $f$ and
rings of the type $S^1\times (0,1)$. Then the cohomology groups
$H_{\cP}^*(Q^3)$ are topological invariants of the system $v$,
i.e., to the topologically equivalent systems there correspond the
isomorphic groups.
\end{th}

       An analogous theorem holds in the multidimensional case as well.
Here the partitioning $\cP$ consists of Liouville tori $T^n$,
$(n-1)$-dimensional tori $T^{n-1}$ and rings of the type
$T^{n-1}\times (0,1)$.

        Since $A^0(\ol{S}(X^{n+1}))=C^\infty(X^{n+1})$, we have
\begin{gather*}
    H_{\cP}^0(X^{n+1})=\ker(d^0:C^\infty(X^{n+1})\to
        A^1(\ol{S}(X^{n+1}))= \\
    =\big\{g\in C^\infty(X^{n+1})\,|\,Xg=0,\;\;
        \forall X\in \ol{S}(X^{n+1}\big\}=
        \Inv_{\ol{S}}C^\infty(X^{n+1}),
\end{gather*}
where $\Inv_{\ol{S}}C^\infty(X^{n+1})$ denotes the set of smooth
functions on $X^{n+1}$ whose restrictions are constant functions
on the elements of the partitioning  $\cP$.

        Denote by $G$ the factor set of the space $X^{n+1}$ with respect
to Liouville tori $T^n$ and the connected components of critical
fibers of the integral $f$. If we introduce the factor topology on
$G$, then $H_{\cP}^0(X^{n+1})$ will coincide with the set of
continuous functions on $G$ whose liftings to $X^{n+1}$ by the
natural projection $X^{n+1}\to G$ are smooth functions.

%        REMARK 3.5.
\begin{rem}
For the four-dimensional case, a complete topological invariant
was introduced in [6]. This is the so-called labelled molecule
consisting of a graph with edges to which are attached rational
numbers from $[0,1)$ or $\infty$.  Note that in the
four-dimensional case the above-mentioned $G$ coincides with the
graph-molecule from [6]. Therefore the zero-dimensional groups
$H_{\cP}^0(X^{n+1})$ already pick up "nonzero" information on the
topological equivalence of integrable Hamiltonian systems.
\end{rem}

%        \S 4.
\section{The Case of a Spherical Pendulum}

        For a spherical pendulum the phase space is a cotangent bundle to
the two-dimensional sphere $T^*S^2$, where $S^2=\{x\in
{\bR}^3:x_1^2+x_2^2+x_3^2=1\}$. Using the Riemannian metric on
$S^2$, we can identify $T^*S^2$ with the tangent bundle $TS^2$ and
define the Hamiltonian of the system by the energy  function
%   (4.1)
\begin{equation}
    E(x,v)=\frac{1}{2}\langle v,v \rangle+x_3,\quad x\in S^2,
        \quad v\in T_xS^2.
\end{equation}

        A kinetic moment with respect to the $x_3$-axis has the form
%   (4.2)
\begin{equation}
    I(x,v)=x_1v_2-v_1x_2.
\end{equation}

        The critical points of the moment mapping $F=(E,I):TS^2\to {\bR}^2$
are the points $x=\pm(0,0,1)$, $v=0$ and $v=\al(-x_2,x_1,0)$,
$1+\al^2x_3=0$, $x_3\neq \pm 1$. The corresponding singular values
of the moment mapping are
%   (4.3)
\begin{equation}
    F=(\pm 1,0) \;\; \text{and} \;\;
        F=\Big(\frac{1}{2}\al^2-\frac{3}{2}\al^{-2},\al-\al^{-3}\Big),
        \;\; \text{where} \;\; |\al|>1.
\end{equation}
 When $\al$ tends to $\pm 1$, the point of the curve defined by
formula (4.3) tends to $(-1,0)$. For $I\neq0$, i.e., for $x\neq
\pm (0,0,1)$, we can introduce the polar coordinates
$$  x_1=\sin\vf\cos\tht, \quad x_2=\sin\vf\sin\tht,\quad x_3=\cos\vf, $$
where $\tht\in [0.2\pi]$, $\vf\in ]0,\pi[$.

        The functions $E$ and $I$ in terms of these coordinates are written as
%   (4.4)
\begin{equation}
    E=\frac{1}{2}\,\dot{\vf}^2+V_I(\vf), \quad I=(\sin^2\vf)\dot{\tht},
\end{equation}
where $V_I(\vf)=\frac{1}{2}(\sin^2\vf)\dot{\tht}^2+\cos\vf$ is the
effective potential, and $(\vf,\tht,\dot{\vf},\dot{\tht})$ define
the local coordinate system on $TS^2$.  The image $F$ is given  by
the relation $I=\al-\al^{-3}$, $E\geq
\frac{1}{2}\al^2-\frac{3}{2}\al^{-2}$, $|\al|\geq 1$. For regular
values of the mapping $F$ we have
$E>\frac{1}{2}\al^2-\frac{3}{2}\al^{-2}$; the point $(0,1)$ should
be discarded. %(see Fig. 1).
%\bigskip
%\input tevd-fig.pic

%\centerline{Fig. 1}
%\smallskip
        Let us consider the function $I$ on the isoenergetic surface
$Q^3=\{E=\frac{1}{2}\}$. If $(\vf,\tht,\dot{\vf},\dot{\tht})$ is a
critical point of the function $I$ on $Q^3$, then at that point
$\grad E$ and $\grad I$ are colinear. Hence we obtain
$\vf_0=\arccos(\frac{1-\sqrt{13}}{6})$,
$\dot{\tht}_1=\sqrt{\frac{\sqrt{13}+1}{2}}$,
$\dot{\tht}_2=-\sqrt{\frac{\sqrt{13}+1}{2}}$. Thus gives us two
isolated critical circles
$S_1=\{(\vf_0,\tht,\dot{\vf}=0,\dot{\tht}_1)| \tht\in [0,2\pi]\}$
and $S_2=\{(\vf_0,\tht,\dot{\vf}=0,\dot{\tht}_2)| \tht\in
[0,2\pi]\}$. Using the method of constant Lagrange multipliers, we
can conclude that $S_1$ is the maximal critical circle, and $S_2$
is the minimal one.

        By Morse's  lemma it follows that the compact isoenergetic surface
corresponding to the values of energy $E$ close to $-1$ is the
three-dimensional
 sphere $S^3$. Therefore $Q^3$ is also diffeomorphic to $S^3$, since
the points $\pm 1$ are critical values of the integral $E$. The
isoenergetic surface $Q^3$ can be represented as two solid tori
$S_1^1\times D^2$ and $S_2^1\times D^2$ glued together along the
boundary tori by means of a diffeomorphism which transforms the
parallel to the meridian, and vice versa (here $D^2$ is a
two-dimensional disk, and the central circles of the  solid tori
$S_1^1\times \{0\}$ and $S_2^1\times \{0\}$ coincide with the
critical circles of the integral $I$); the tori $S_i^1\times S_r$,
$r\in (0,1]$, $i=1,2$, where $S_r=\{x\in D^2\,|\,\|x\|=r\}$, are
the usual Liouville tori of the Hamiltonian system $v=s\grad E$.
In that case, using the notation from [6], the labelled molecule
has the form
$$  A \bullet\!\!\!\us{r=0}{\text{---\!\!---\!\!---\!\!---}}\!\!\!\bullet A, $$
and the factor set $G$ coincides with the segment.

        Now consider the function $I$ on the isoenergetic surface
$Q^3=\{E=\al>1\}$. A direct calculation shows that, as in the
preceding case, for the integral $I$ we have two isolated critical
circles $S_1=\{(\vf_0,\tht,0,\dot{\tht}_1)\,|\,\tht\in [0,2\pi]\}$
and $S_2=\{(\vf_0,\tht,0,\dot{\tht}_2)\,|\,\tht\in [0,2\pi]\}$,
where $\vf_0=\arccos\frac{\al-\sqrt{\al^2+12}}{3}$ and
$\dot{\tht}_1=-\sqrt{\frac{3}{\sqrt{\al^2+12}-\al}}$,
$\dot{\tht}_2=\sqrt{\frac{3}{\sqrt{\al^2+12}-\al}}$.

        In terms of the coordinates
$(\vf,\tht,\dot{\vf},\dot{\tht})$ the Euclidean metric is written
as $ds^2=d\vf^2+\sin^2\vf d\tht^2$ so that the the equation
$E=\al$ takes the form
$$  \|\xi\|^2+2\cos\vf=2\al,    $$
where $\xi=(\dot{\vf},\dot{\tht})$ is the tangent vector at the
point $(\vf,\tht)\in S^2$, and $\|\xi\|$ denotes the norm of
$\xi$. Hence it follows that the norm of the vector $\xi$ is a
function of $\vf$,and $\|\xi\|\neq 0$. Therefore $Q^3$ is
diffeomorphic to the space  $T_1S^2$, where
$T_1S^2=\{(a,\xi)\,|\,x\in S^2$, $\|\xi\|=1\}$. As is well-known,
$T_1S^2$ is diffeomorphic to the three-dimensional projective
space ${{\bR} P}^3$ so that $Q^3\approx {{\bR} P}^3$. As above,
$Q^3$ can be represented as two solid tori $S_1^1\times D^2$ and
$S_2^1\times D^2$ glued together along the boundary tori by means
of the diffeomorphism $h$ whose corresponding induced mapping
$h_*:H^1(T^2;{\bZ})\to H^1(T^2;{\bZ})$ is given by the matrix
$\begin{pmatrix} 1 & 2 \\ 1 & 1 \end{pmatrix}$. Thus for the
integrable Hamiltonian system $v=\sgrad E$ on $Q^3$ the labelled
molecule has the form
$$  A \bullet \!\!\!\us{r=\frac{1}{2}}{\text{---\!\!---\!\!---\!\!---}}\!\!\!\bullet A, $$
and the factor set again coincides with the segment. This example
shows that vertical cohomologies are not a complete topological
invariantof integrable Hamiltonian systems.

\section*{Acknowledgements}

The author expresses his gratitude to Dr. Z. Giunashvili for
helpful discussions.

%This research was partly carried out under Grant No.
%RJV200 from  the Intenational Science Foundation.

\vskip+0.8cm

\centerline{\sc References}
\vskip+0.4cm

1. Jos\'{e} M. Figueroa-o'Farrill, A topological characterization of classical
    BRST cohomology. {\em Commun. Math. Phys.}
    {\bf 127}(1990), 181--186.

2. M. Henneaux and C. Teitelboim,. BRST cohomology in classical mechanics.
    {\em Commun. Math. Phys.} {\bf 115}(1988), 213--230.

3. J. Fisch, M. Henneaux, J. Stasheff, and C. Teitelboim, Existence, uniqueness,
    and cohomology of the classical BRST charge with ghosts of ghosts.
    {\em Preprint,} 1988.

4. I. Tamura, Topology foliations. {\em Japan: Iwanami Shoten,} 1976.

5. A. T. Fomenko, The topology of surfaces of constant energy in integrable
    hamiltonian systems, and obstructions to integrability.
    {\em Math. USSR Izvestiya} {\bf  29}(1987), No. 3, 629--658.

6. A. V. Bolsinov, S. V. Matveev, and A. T. Fomenko, Topological classification of
    integrable Hamiltonian systems with two degrees of freedom.
    A list of systems of small complexity.
    (Translation from Russian) {\em  Math. Surveys} {\bf 45}(1990), No. 2,
    49--77.

\

\

%\centerline{(Received 1.03.1996)}

\

Author's address:

Faculty of Mechanics and Mathematics

I. Javakishvili Tbilisi State University

2, University St., Tbilisi 380043  Georgia

 E-mail address:
tevzaz@@ti.net.ge
\end{document}